\documentclass{custom}
 \usepackage{amsmath, amssymb, bm, mathtools}

 \def\TS{\textstyle} \def\d{{\rm d}}  
 \def\VEC#1{{\pmb{#1}}} \def\MAT#1{{\pmb{#1}}} \def\IE{{\it i.e.}} 
 \def\EG{{\it e.g.}}  
 \def\FIRST{{$1^{\textrm{st}}$}} \def\SECOND{{$2^{\textrm{nd}}$}} 
 \def\THIRD{{$3^{\textrm{rd}}$}} \def\NTH#1{{${#1}^{\textrm{th}}$}} 
 \def\DOTS{{...}} \def\STRUT{\vphantom{|_|^|}} 
 
 \def\AIJ{ a_{\mkern1mu i \mkern-1mu j} } 
 \def\BJ{ b_{\mkern-1.25mu j} }
 \def\CJ{ c_{\mkern-1.25mu j} }
 \def\FJ{ \VEC{F}_{\mkern-3.5mu j} }
 \def\SJ{ \beta_{\mkern-0.25mu j} }
 \def\SKJ{ \beta_{k \mkern-1mu j} } 
 \def\YJ{ Y_{\mkern-1mu j} }

\begin{document}

  \title{Two continuous $\textbf{(4, 5)}$ pairs of explicit 
$\textbf{9}$-stage Runge--Kutta methods}
  \shorttitle{Two continuous $(4, 5)$ pairs of Runge--Kutta methods}
  \author{{\sc Misha Stepanov\thanks{Email: stepanov@math.arizona.edu}} 
\\[2pt] Department of Mathematics \;and\; Program in Applied 
Mathematics,{}\\{}University of Arizona, Tucson, AZ 85721, USA}
  \shortauthorlist{M.~Stepanov}

\maketitle

\begin{abstract} {An $11$-dimensional family of embedded $(4, 5)$ pairs 
of explicit $9$-stage Runge--Kutta methods with an interpolant of order 
$5$ is derived. Two optimized for efficiency pairs are presented.} 
{embedded pair of Runge--Kutta methods, continuous formula, interpolant} 
\end{abstract}

\section{Introduction}

Runge--Kutta methods (see, \EG, \citep[][sec.~23 and ch.~3]{But08}, 
\citep[][ch.~II]{HNW08}, \citep[][ch.~4]{AsPe98}, \citep[][ch.~3]{Ise09}) 
are widely and successfully used to solve Ordinary Differential Equations 
(ODEs) numerically for over a century \citep[][]{BuWa96}. Being applied 
to a system $\d \mkern0.25mu \VEC{x} / \d \mkern0.25mu t = \VEC{f}(t, \, 
\VEC{x})$, in order to propagate by the step size $h$ and update the 
position, $\VEC{x}(t) \mapsto \VEC{x}(t + h)$, an \mbox{$s$-stage} 
explicit Runge-Kutta method (which is determined by the coefficients 
$\AIJ$, weights $b_{\mkern-1mu j}$, and nodes $c_{\mkern0.5mu i}$) would 
compute intermediate vectors $\VEC{F}_{\mkern-1mu 1}$, 
$\VEC{X}_{\mkern-1mu 2}$, $\VEC{F}_{\mkern-1mu 2}$, $\VEC{X}_{\mkern-1mu 
3}$, \DOTS, $\VEC{X}_{\mkern-1mu s}$, $\VEC{F}_{\mkern-3mu s}$, and then 
$\VEC{x}(t + h)$:
 \begin{gather*} \VEC{X}_{\mkern-1mu i} = \VEC{x}(t) + h \! \sum_{j = 
1}^{i - 1} \AIJ \mkern2mu \FJ, \qquad \VEC{F}_{\mkern-1.5mu i} = \VEC{f} 
\bigl( t + c_{\mkern0.5mu i} \mkern1mu h, \mkern2mu \VEC{X}_{\mkern-1mu 
i} \bigr), \qquad \VEC{x}(t + h) = \VEC{x}(t) + h \! \sum_{j = 1}^s \BJ 
\mkern2mu \FJ \end{gather*} In the limit $h \to 0$ all the vectors 
$\VEC{F}_{\mkern-1.5mu i}$, where $1 \le i \le s$, are the same, so it 
is natural and will be assumed that $\sum_{\mkern-1mu j = 
1}^{\mkern0.5mu i - 1} \AIJ = c_{\mkern0.5mu i}$. For $i = 1$ the sum 
over $\mkern-1mu j$ is empty, so $c_1 = 0$, $\VEC{X}_{\mkern-1mu 1} = 
\VEC{x}(t)$, and $\VEC{F}_{\mkern-1mu 1} = \smash{\VEC{f} \bigl( t, 
\mkern1mu \VEC{x}(t) \bigr)}$.

To obtain an accurate solution with less effort, various adaptive step 
size strategies were developed (see, \EG, \citep[][secs.~271 and 
33]{But08}, \citep[][sec.~II.4]{HNW08}, \citep[][sec.~4.5]{AsPe98}, 
\citep[][ch.~6]{Ise09}). Typically a system of ODEs is solved in two 
different ways, and the step size is chosen so that the two solutions 
are sufficiently close. A computationally efficient procedure is to have 
two Runge--Kutta methods with different weights, but the same nodes and 
coefficients. The vectors $\VEC{F}_1$, $\VEC{F}_2$, \DOTS, 
$\VEC{F}_{\mkern-3mu s}$ are computed only once, and then are used in 
both methods, the latter are said to form an embedded pair. Two well 
known examples of such pairs are \citep[][tab.~III]{Feh69}, 
\citep[][tab.~1]{Feh70} and \citep[][tab.~2]{DoPr80}.

There is no $5$-stage explicit Runge--Kutta \NTH{5} order method 
\citep{But64}. The Fehlberg pair has $6$ stages. The Dormand--Prince 
pair uses $7$ stages, but has the so-called First Same As Last (FSAL) 
property \citep[p.~17]{Feh69}, \citep{DoPr78}: the vector 
$\VEC{F}_{\mkern-1mu 1}$ at the current step is equal to the already 
computed $\VEC{F}_{\mkern-2mu u}$ at the stage $1 < u \le s$ of the 
previous step. An FSAL method requires $(s - 1)$ evaluations of the 
r.h.s.~function $\VEC{f}$ per step, with the exception of the \FIRST{} 
step.

  \begin{table}[t] \centerline{\begin{tabular}{ccccc} 
    {\hskip0.05in{order}\hskip0.05in} & {\hskip0.45in{$t$}\hskip0.45in} 
& {\hskip0.05in{$\gamma(t)$}\hskip0.05in} & 
{\hskip0.05in{$\sigma(t)$}\hskip0.05in} & 
{\hskip0.45in{$\VEC{\Phi}(t)$}\hskip0.45in} \\
    \hline
    \FIRST{} & \begin{picture}(0,12)(0,2.5) \thicklines
  \put(0,6){\circle*{3}}
\end{picture} & $1$ & $1$ & $\VEC{1} \STRUT$ \\
    \hline
    \SECOND{} & \begin{picture}(12,12)(0,3) \thicklines
  \put(0,6){\line(1,0){12}}
  \put(0,6){\circle*{3}}
  \put(12,6){\circle*{3}}
\end{picture} & $2$ & $1$ & $\VEC{c} \STRUT$ \\
    \hline
    \THIRD{} & \begin{picture}(12,12)(0,3) \thicklines
  \put(0,6){\line(4,-1){12}}
  \put(0,6){\line(4,1){12}}
  \put(0,6){\circle*{3}}
  \put(12,3){\circle*{3}}
  \put(12,9){\circle*{3}}
\end{picture} & $3$ & $2$ & $\VEC{c}^2$ \\
    \raisebox{4pt}{$\downarrow$} &
\begin{picture}(24,12)(0,3) \thicklines
  \put(0,6){\line(1,0){24}}
  \put(0,6){\circle*{3}}
  \put(12,6){\circle*{3}}
  \put(24,6){\circle*{3}}
\end{picture} & $6$ & $1$ & $\MAT{A} \mkern1mu \VEC{c} = \VEC{c}^2 
\mkern-1mu / 2 + \VEC{q}_1 \STRUT$ \\
    \hline
    \NTH{4} & \begin{picture}(12,12)(0,3) \thicklines
  \put(0,6){\line(1,0){12}}
  \put(0,6){\line(3,-1){12}}
  \put(0,6){\line(3,1){12}}
  \put(0,6){\circle*{3}}
  \put(12,2){\circle*{3}}
  \put(12,6){\circle*{3}}
  \put(12,10){\circle*{3}}
\end{picture} & $4$ & $6$ & $\VEC{c}^3$ \\
  $\downarrow$ & \begin{picture}(24,12)(0,3) \thicklines
  \put(0,6){\line(4,-1){12}}
  \put(0,6){\line(4,1){12}}
  \put(12,3){\line(1,0){12}}
  \put(0,6){\circle*{3}}
  \put(12,3){\circle*{3}}
  \put(12,9){\circle*{3}}
  \put(24,3){\circle*{3}}
\end{picture} & $8$ & $1$ & $(\MAT{A} \mkern1mu \VEC{c}) \mkern1mu 
\VEC{c} = \VEC{c}^3 \mkern-1mu / 2 + c_2 \mkern1mu \VEC{q}_1$ \\
    & \begin{picture}(24,12)(0,3) \thicklines
  \put(0,6){\line(1,0){12}}
  \put(12,6){\line(4,-1){12}}
  \put(12,6){\line(4,1){12}}
  \put(0,6){\circle*{3}}
  \put(12,6){\circle*{3}}
  \put(24,3){\circle*{3}}
  \put(24,9){\circle*{3}}  
\end{picture} & $12$ & $2$ & $\MAT{A} \mkern1mu \VEC{c}^2 = \VEC{c}^3 
\mkern-1mu / 3 + \VEC{q}_2$ \\
    & \begin{picture}(36,12)(0,3) \thicklines
  \put(0,6){\line(1,0){36}}
  \put(0,6){\circle*{3}}
  \put(12,6){\circle*{3}}
  \put(24,6){\circle*{3}}
  \put(36,6){\circle*{3}}
\end{picture} & $24$ & $1$ & $\MAT{A}^2 \mkern1mu \VEC{c} = \MAT{A} 
\mkern1mu \VEC{c}^2 / 2 + \MAT{A} \mkern1mu \VEC{q}_1 \STRUT$ \\
    \hline
    \NTH{5} & \begin{picture}(12,12)(0,3) \thicklines
  \put(0,6){\line(2,-1){12}}
  \put(0,6){\line(6,-1){12}}
  \put(0,6){\line(6,1){12}}
  \put(0,6){\line(2,1){12}}
  \put(0,6){\circle*{3}}
  \put(12,0){\circle*{3}}
  \put(12,4){\circle*{3}}
  \put(12,8){\circle*{3}}
  \put(12,12){\circle*{3}}
\end{picture} & $5$ & $24$ & $\VEC{c}^4$ \\
    $\downarrow$ & \begin{picture}(24,12)(0,3) \thicklines
  \put(0,6){\line(1,0){12}}
  \put(0,6){\line(3,-1){12}}
  \put(0,6){\line(3,1){12}}
  \put(12,2){\line(1,0){12}}
  \put(0,6){\circle*{3}}
  \put(12,2){\circle*{3}}
  \put(12,6){\circle*{3}}
  \put(12,10){\circle*{3}}
  \put(24,2){\circle*{3}}
\end{picture} & $10$ & $2$ & $(\MAT{A} \mkern1mu \VEC{c}) \, \VEC{c}^2 = 
\VEC{c}^4 \mkern-1mu / 2 + c_2^2 \mkern1mu \VEC{q}_1$ \\
    & \begin{picture}(24,12)(0,3) \thicklines
  \put(0,7.5){\line(4,-1){12}}
  \put(0,7.5){\line(4,1){12}}
  \put(12,4.5){\line(4,-1){12}}
  \put(12,4.5){\line(4,1){12}}
  \put(0,7.5){\circle*{3}}
  \put(12,4.5){\circle*{3}}
  \put(12,10.5){\circle*{3}}
  \put(24,1.5){\circle*{3}}
  \put(24,7.5){\circle*{3}}
\end{picture} & $15$ & $2$ & $( \MAT{A} \mkern1mu \VEC{c}^2 ) \mkern1mu 
\VEC{c} = \VEC{c}^4 \mkern-1mu / 3 + \VEC{q}_2 \mkern1mu \VEC{c}$ \\
    & \begin{picture}(36,12)(0,3) \thicklines
  \put(0,6){\line(4,-1){12}}
  \put(0,6){\line(4,1){12}}
  \put(12,3){\line(1,0){24}}
  \put(0,6){\circle*{3}}
  \put(12,3){\circle*{3}}
  \put(12,9){\circle*{3}}
  \put(24,3){\circle*{3}}
  \put(36,3){\circle*{3}}
\end{picture} & $30$ & $1$ & $(\MAT{A}^2 \mkern1mu \VEC{c}) \mkern1mu 
\VEC{c} = (\MAT{A} \mkern1mu \VEC{c}^2) \, \VEC{c} / 2 + (\MAT{A} 
\mkern1mu \VEC{q}_1 \mkern-1.5mu ) \mkern1mu \VEC{c}$ \\
    & \begin{picture}(24,12)(0,3) \thicklines
  \put(0,6){\line(4,-1){12}}
  \put(0,6){\line(4,1){12}}
  \put(12,3){\line(1,0){12}}
  \put(12,9){\line(1,0){12}}
  \put(0,6){\circle*{3}}
  \put(12,3){\circle*{3}}
  \put(12,9){\circle*{3}}
  \put(24,3){\circle*{3}}
  \put(24,9){\circle*{3}}
\end{picture} & $20$ & $2$ & $(\MAT{A} \mkern1mu \VEC{c})^2 = \VEC{c}^4 
\mkern-1mu / 4 + c_2^2 \mkern1mu \VEC{q}_1 / 2$ \\
    & \begin{picture}(24,12)(0,3) \thicklines
  \put(0,6){\line(1,0){12}}
  \put(12,6){\line(1,0){12}}
  \put(12,6){\line(3,-1){12}}
  \put(12,6){\line(3,1){12}}
  \put(0,6){\circle*{3}}
  \put(12,6){\circle*{3}}
  \put(24,2){\circle*{3}}
  \put(24,6){\circle*{3}}
  \put(24,10){\circle*{3}}
\end{picture} & $20$ & $6$ & $\MAT{A} \mkern1mu \VEC{c}^3 = \VEC{c}^4 
\mkern-1mu / 4 + \VEC{q}_3$ \\
    & \begin{picture}(36,12)(0,3) \thicklines
  \put(0,6){\line(1,0){12}}
  \put(12,6){\line(4,-1){12}}
  \put(12,6){\line(4,1){12}}
  \put(24,3){\line(1,0){12}}
  \put(0,6){\circle*{3}}
  \put(12,6){\circle*{3}}
  \put(24,3){\circle*{3}}
  \put(24,9){\circle*{3}}
  \put(36,3){\circle*{3}}
\end{picture} & $40$ & $1$ & $\MAT{A} \bigl( \mkern-1mu (\MAT{A} 
\mkern1mu \VEC{c}) \mkern1mu \VEC{c} \mkern-0.5mu \bigr) = \MAT{A} 
\mkern1mu \VEC{c}^3 \mkern-1mu / 2 + c_2 \mkern1mu \MAT{A} \mkern1mu 
\VEC{q}_1$ \\
    & \begin{picture}(36,12)(0,3) \thicklines
  \put(0,6){\line(1,0){24}}
  \put(24,6){\line(4,-1){12}}
  \put(24,6){\line(4,1){12}}
  \put(0,6){\circle*{3}}
  \put(12,6){\circle*{3}}
  \put(24,6){\circle*{3}}
  \put(36,3){\circle*{3}}
  \put(36,9){\circle*{3}}
\end{picture} & $60$ & $2$ & $\MAT{A}^2 \mkern1mu \VEC{c}^2 = \MAT{A} 
\mkern1mu \VEC{c}^3 \mkern-1mu / 3 + \MAT{A} \mkern1mu \VEC{q}_2$ \\
    & \begin{picture}(48,12)(0,3) \thicklines
  \put(0,6){\line(1,0){48}}
  \put(0,6){\circle*{3}}
  \put(12,6){\circle*{3}}
  \put(24,6){\circle*{3}}
  \put(36,6){\circle*{3}}
  \put(48,6){\circle*{3}}
\end{picture} & $120$ & $1$ & $\MAT{A}^3 \mkern1mu \VEC{c} = \MAT{A}^2 
\mkern1mu \VEC{c}^2 / 2 + \MAT{A}^2 \mkern1mu \VEC{q}_1$ \end{tabular}} 
\caption{Order conditions $\VEC{b} \mkern1mu \VEC{\Phi}(t) = 1 / 
\gamma(t)$ for rooted trees $t$ with up to $5$ vertices. It is assumed 
that $(\MAT{A} \mkern1mu \VEC{c})_{\mkern0.5mu i} = c_{\mkern0.5mu i}^2 
/ 2$ for all $i \ne 2$.} \label{order_conditions} \end{table}

Continuous formulas or interpolants (see, \EG, \citep[][]{Hor83}, 
\cite[][]{Sar84}, \citep[][sec.~272]{But08}, \citep[][sec.~II.6]{HNW08}) 
provide an inexpensive (\IE, with only a few if any additional 
evaluations of the r.h.s.) way to estimate the solution at anywhere 
within the integration interval. Without altering the strategy of step 
size choice, this can be used in applications that require values of the 
solution $\VEC{x}(t)$ at specific points $t_1$, $t_2$, \DOTS{} (dense 
output) or the place (time $t$ and/or position $\VEC{x}(t)$) where the 
solution crosses a hypersurface $\smash{g \bigl( t, \mkern1mu \VEC{x}(t) 
\bigr) = 0}$ (event location). The continuous approximation to the 
solution in the interval $[ \mkern1mu t,\, t + h \mkern1mu ]$ is 
typically of the form \begin{gather*} \VEC{x}(t + \theta h) = \VEC{x}(t) 
+ h \! \sum_{j = 1}^s \SJ(\theta) \mkern2mu \VEC{F}_{\mkern-3.5mu j} 
\end{gather*} where the interpolant functions $\SJ(\theta) = \sum_k \SKJ 
\theta^k$ are polynomials. For the approximation over several step 
intervals to be continuously differentiable a method should have the FSAL 
property, with the following conditions on the behavior of the row vector 
$\VEC{\beta}(\theta) = \bigl[ \SJ(\theta) \bigr]$ at $\theta = 0$ and 
$\theta = 1$:
 \begin{align*} c_1 &= 0, & \VEC{a}_{1*} &= \VEC{0} = \VEC{\beta}(0), & 
\VEC{X}_{\mkern-1mu 1} &= \VEC{x}(t), & \VEC{F}_{\mkern-1mu 1} &= \VEC{f} 
\bigl( t, \, \VEC{x}(t) \bigr), & \frac{\d \SJ(\theta)}{\d \theta} 
\Bigr|_{\theta = 0} &= \left\{ \begin{array}{ll} 1, & j = 1 \\ 0, & j \ne 
1 \end{array} \right. \\ c_u &= 1, & \VEC{a}_{u*} &= \VEC{b} = 
\VEC{\beta}(1), & \VEC{X}_{\mkern-1mu u} &= \VEC{x}(t + h), & 
\VEC{F}_{\mkern-2mu u} &= \VEC{f} \bigl( t + h, \, \VEC{x}(t + h) \bigr), 
& \frac{\d \SJ(\theta)}{\d \theta} \Bigr|_{\theta = 1} &= \left\{ 
\begin{array}{ll} 1, & j = u \\ 0, & j \ne u \end{array} \right.
 \end{align*} Here $\VEC{a}_{i*} = \bigl[ \AIJ \bigr]$ and $\VEC{b} = 
\bigl[ \BJ \bigr]$ are the row vectors of coefficients and weights, 
respectively. Below all interpolants are assumed to be continuously 
differentiable.

The product of column vectors $\VEC{x} \mkern2mu \VEC{y}$ is to be 
understood element-wise: $(\VEC{x} \mkern2mu \VEC{y})_i = x_i \mkern1mu 
y_i$. Let $\VEC{1}$ be the $s$-dimensional column vector with all 
components being equal to $1$; $\VEC{c} = \bigl[ c_i \bigr]$ be the nodes 
vector; and $\MAT{A} = \bigl[ \AIJ \bigr]$ be the $s \times s$ matrix 
with $\AIJ$ as its matrix element in the $i^{\mkern1mu \textrm{th}}$ row 
and $j^{\mkern1mu \textrm{th}}$ column (for an explicit method $\AIJ = 0$ 
if $i \le j$). Let $\VEC{q}_n = \MAT{A} \mkern1mu \VEC{c}^n - \frac{1}{n 
+ 1} \mkern1mu \VEC{c}^{n + 1}$. The condition $\sum_j \AIJ = c_{i}$, 
$\MAT{A} \VEC{1} = \VEC{c}$, or $\VEC{q}_0 = \VEC{0}$ is assumed. The 
following quantities will be used for the estimation of the local error: 
\begin{align*}
   T_p^2(\VEC{x}, \theta) = 
\sum_{\mathclap{\text{rooted~trees~}t\text{~of~order~}p}} \tau^2(t, 
\mkern2mu \VEC{x}, \theta), \qquad
  \tau(t, \mkern2mu \VEC{x}, \theta) = \frac{1}{\sigma(t)} \biggl( 
\VEC{x} \mkern1.5mu \VEC{\Phi}(t) - \frac{\theta^p}{\gamma(t)} \biggr)
 \end{align*} Here $\sigma(t)$ is the order of the symmetry group of the 
tree $t$ (see, \EG, \cite[p.~140]{But08}). Whenever the argument 
$\VEC{x}$ or $\theta$ is omitted, its value is meant to be equal to 
$\VEC{\beta}(\theta)$ and $1$, respectively. Note that $\VEC{\beta}(1) = 
\VEC{b}$. An interpolant of order $p$ should satisfy the conditions 
$\tau(t, \theta) \equiv 0$ for all rooted trees $t$ with up to $p$ 
vertices. For $p = 5$, with the assumption that the intermediate 
positions $\VEC{X}_{\mkern-1mu i}$ are at least \SECOND{} order accurate 
for all $i > 2$, these conditions are listed in 
Table~\ref{order_conditions}, see also \citep[][sec.~31]{But08}, 
\citep[][sec.~II.2]{HNW08}, \citep[][tab.~1]{DoPr80}.

There are certain properties one would expect from a practical 
Runge--Kutta method (see, \EG, a list in \citep[][p.~785]{Ver78}). The 
natural or desirable behavior of the interpolant function $\SJ(\theta)$, 
where $1 \le j \le s$, is that it has a notably positive slope around 
$\theta \approx \CJ$, while it is hardly changing anywhere else. 
Deviations from this behavior can be divided into two categories: {\em 
non-positivity}, when an interpolant function's derivative is negative; 
and {\em non-locality}, when an interpolant function has substantial 
slope (positive or negative) far from the position of the corresponding 
node.

\begin{figure} \centerline{\includegraphics{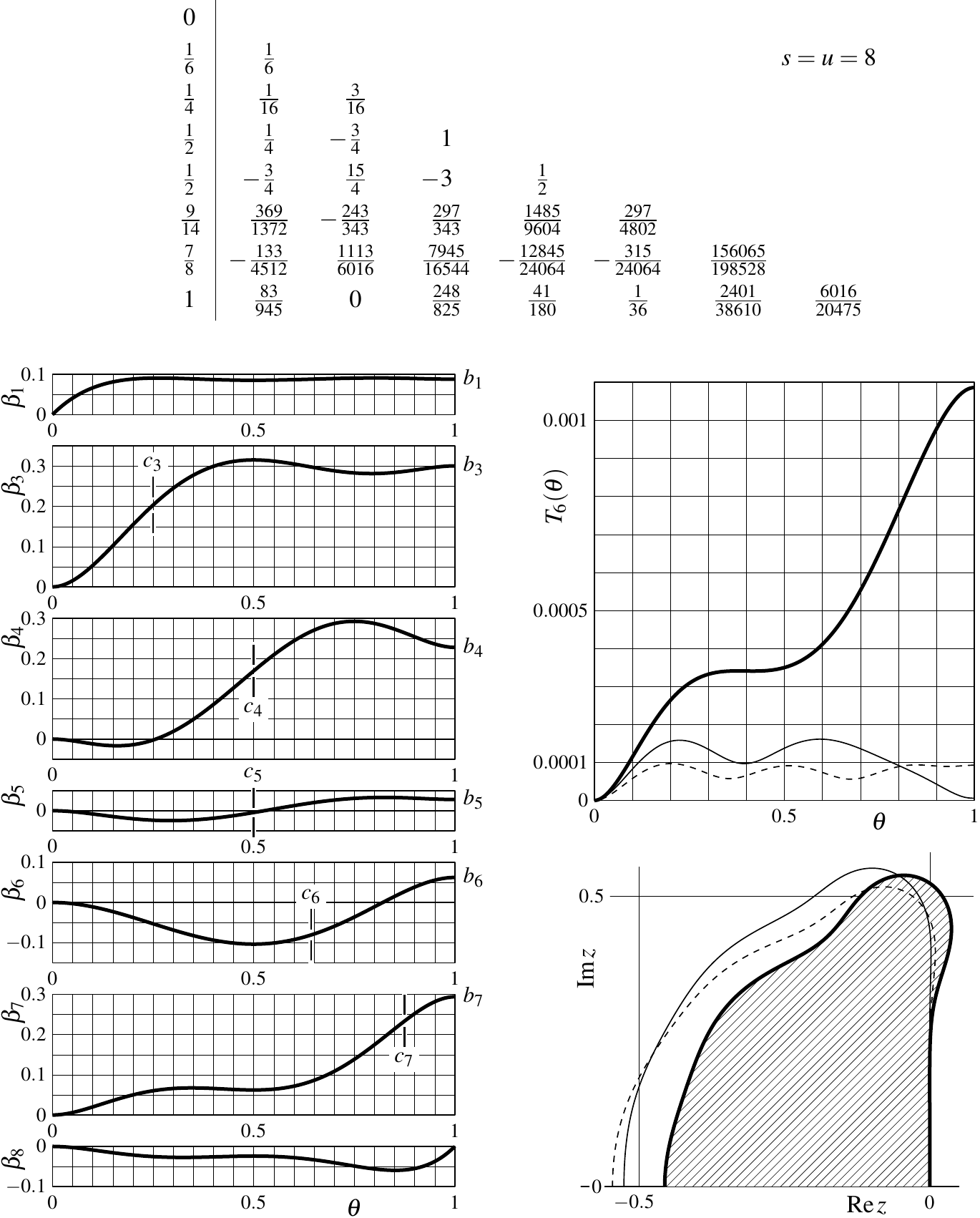}} 
\caption{The continuous $(4, 5)$ pair \citep[][fig.~3]{OwZe92}: the 
Butcher tableau (upper panel), the weights vector $\VEC{b}$ repeats the 
last row of $\MAT{A}$ and is not shown; the interpolant functions 
$\SJ(\theta)$, with $j = 1$, $3$, $4$, \DOTS, $8$ (the panels on the 
left, the function $\beta_{\mkern1.5mu 2}(\theta) \equiv 0$ is not 
shown); the local error $T_6(\theta)$ (middle right panel); and the 
region of absolute stability (lower right panel). (The thin dashed and 
solid lines correspond to the pairs on Figure~\ref{pair1} and 
Figure~\ref{pair2}, respectively. The three regions of absolute stability 
are scaled for equal cost, \IE, the regions where $\bigl| \mkern1mu R 
\bigl( \mkern-1mu (s - 1) \mkern1mu z \bigr) \bigr| \le 1$ are depicted, 
here $R(z)$ is the stability function and $s$ is the number of stages.)} 
\label{OwZe} \end{figure}%
\hskip\parindent An explicit Runge--Kutta method with continuously 
differentiable interpolant of order $5$ has at least $8$ stages 
\citep[][sec.~3.3]{OwZe91}. Such methods were completely classified in 
\citep[][]{VeZe95}. In \citep[][]{OwZe92} a $5$-dimensional family of 
continuous $(4, 5)$ pairs of $8$-stage Runge-Kutta methods was 
constructed, and an optimized for efficiency \cite[][fig.~3]{OwZe92} pair 
was suggested, which is also shown in Figure~\ref{OwZe}. All the pairs in 
this $5$-dimensional family satisfy $2 c_3 = c_4 = c_5$. This is an 
indication that the family lacks sufficient flexibility or is stressed by 
numerous imposed conditions. Looking at the interpolant functions in 
Figure~\ref{OwZe}, $\beta_4(\theta)$ goes down for $\theta > 0.8$, 
$\beta_6(\theta)$ goes down for $\theta < 0.5$, and both 
$\beta_7(\theta)$ and $\beta_8(\theta)$ express non-locality for $\theta 
< 0.3$.

Any Runge--Kutta method could be equipped with an interpolant by adding, 
if needed, additional stages \citep[][]{BOOTSTRAP}, \citep[][]{Ver93}. 
Several interpolants were constructed for the \citep[][tab.~2]{DoPr80} 
pair, see, \EG, \citep[][p.~149]{Sha86} and \citep[][]{CMR90}. With $u = 
7$ both interpolants use $s = 9$ stages (see also 
\citep[][corollary~2.13]{OwZe91}). The second interpolant has somewhat 
smaller local error, see \citep[][fig.~1]{CMR90}.

Increasing the number of stages (and thus the amount of computation per 
step) provides additional flexibility in choosing the nodes, 
coefficients, and weights, which may be exploited to construct viable 
pairs that produce an accurate solution in fewer steps. In 
\citep[][sec.~3.1]{ShSm93} and \citep[][]{BoSh96} non-FSAL embedded $(4, 
5)$ pairs of $7$-stage Runge--Kutta methods were suggested. The objective 
in the construction of the latter pair was an improvement of the 
\citep[][tab.~2]{DoPr80} pair (see \citep[][p.~19]{BoSh96}). The pair was 
also equipped with an interpolant of order $5$. The minimal number of 
stages would be $9$, but in \citep[][]{BoSh96} the suggested interpolant 
(with the local error, to the leading order, being a problem-independent 
function of the local error at the end of the step) is using $11$ stages.

\begin{figure} \includegraphics{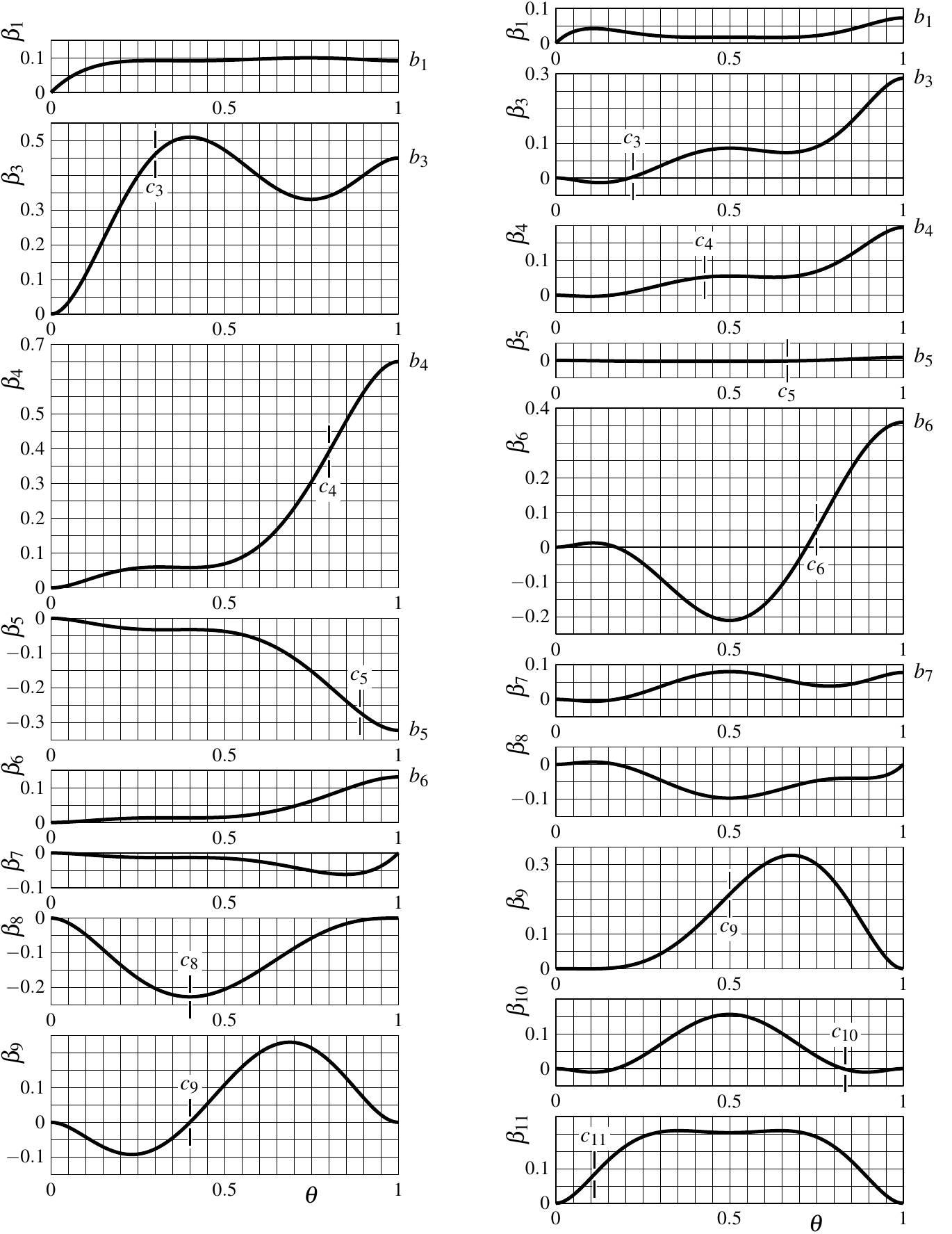} \caption{The 
interpolant functions $\SJ(\theta)$, with $j = 1$, $3$, $4$, \DOTS, $9$ 
\citep[][]{CMR90} for the embedded $(4, 5)$ pair \citep[][tab.~2]{DoPr80} 
(panels on the left, $s = 9$, $u = 7$) and functions $\SJ(\theta)$, with 
$j = 1$, $3$, $4$, \DOTS, $11$, for the \citep[][]{BoSh96} pair (panels 
on the right, $s = 11$, $u = 8$).} \label{DP_BS} \end{figure}%
\hskip\parindent The interpolant \citep[][]{CMR90} for the 
\citep[][tab.~2]{DoPr80} pair and the interpolant for the 
\citep[][]{BoSh96} pair are depicted in Figure~\ref{DP_BS}. In the case 
of Dormand--Prince pair, $b_5 < 0$, plus the \FIRST{} order condition 
$\sum_j \SJ(\theta) \equiv \theta$ is obtained through the cancellation 
of wiggles in $\beta_3$, $\beta_8$, $\beta_9$, and of slopes in $\beta_4$ 
and $\beta_5$ for $\theta > 0.6$. In the case of Bogacki--Shampine pair, 
there is some non-locality in $\beta_1$, $\beta_3$, $\beta_4$, and 
$\beta_6$, which is compensated by $\beta_9$, $\beta_{10}$, and 
$\beta_{11}$. If an interpolant is obtained by adding stages to an 
already formed pair, the interpolant function for an added stage has a 
negative slope somewhere, as it should have zero values (and derivatives) 
at $\theta = 0$ and $\theta = 1$.

In this work embedded $(4, 5)$ pairs are constructed, like the 
\cite[][fig.~3]{OwZe92} pair, so that they have an interpolant right 
away, although not the minimal number of stages is used. The family of 
such pairs is constructed in Section~\ref{eleven}. How the values of the 
free parameters are chosen is discussed in Section~\ref{parameters}. The 
performance of built pairs is demonstrated in Section~\ref{tests}.

\section{A family of continuous $\pmb{(4, 5)}$ pairs} \label{eleven}

There are two different ways to generate an embedded pair with an 
interpolant: to construct a pair with no interpolant, and then add one or 
several stages in order to build one; or to design both a pair and an 
interpolant at once. The latter approach is used here. As for continuous 
$(4, 5)$ pairs $8$ stages do not provide enough flexibility 
\citep[][]{OwZe92}, here $s = 9$ stages are used. For the update of 
position, $\VEC{x}(t) \mapsto \VEC{x}(t + h)$, to be as accurate as 
possible, the FSAL stage is the last one: $u = 9$.

To increase the similarity with collocation methods (see, \EG, 
\citep[][p.~211]{HNW08}, \citep[][sec.~4.7.1]{AsPe98}, 
\citep[][sec.~3.4]{Ise09}), the Dominant Stage Order (DSO) (see, \EG, 
\citep[][eq.~(5)]{Ver10}) is chosen to be equal to $3$. This goes against 
the observation \citep[][p.~386]{Ver10} that most efficient for 
computation pairs of order $p$ have the DSO being equal to $(p - 4)$ or 
$(p - 3)$. Increasing the DSO makes order conditions more redundant, and 
may not reduce richness or flexibility of the set of pairs much.

The parameters of the $11$-dimensional family of continuous $(4, 5)$ 
pairs described below are $c_2$, $c_4$, $c_5$, $c_6$, $c_7$, $c_8$, 
$a_{65}$, $a_{75}$, $a_{76}$, $a_{86}$, and $a_{87}$. Their values are 
arbitrary, except for some degenerate cases, \EG, $c_5 = c_4$ for which 
the matrix $\MAT{A}$ ends up being infinite. Other nodes $c_1$, $c_3$, 
$c_9$, and the first $8$ rows of $\MAT{A}$ are expressed through the $11$ 
parameters as follows: \begin{align*}
   c_1 &= 0, \qquad c_3 = 2 c_4 / 3, \qquad c_9 = 1 \\
   h_{\mkern1mu i \mkern-1mu j} &= \AIJ \mkern1mu \CJ (\CJ - c_4) 
\mkern3mu \bigg/ \mkern18mu \prod^{k \ne i}_{\mathclap{\mkern12mu k \in 
\{1, 4, 5, 6, 7, 8 \}}} \mkern18mu (c_i - c_k), \qquad
   \YJ = 3 - 5 c_4 - 5 \CJ + 10 c_4 \CJ \\
   Z_m &= 12 - 15 c_4 - 15 c_5 - 15 c_m + 20 c_4 c_5 + 20 c_4 c_m + 20 
c_5 c_m - 30 c_4 c_5 c_m \quad \smash{\left. \vphantom{\rule{1pt}{42pt}} 
\right\}} \mkern9mu \mbox{\small\smash{\begin{minipage}{0.8in}
determines $a_{85} \STRUT$\\
ensures $b_9 = 0 \STRUT$
\end{minipage}}} \\
   & \sum_{i \mkern-1mu j \mkern2mu \in \mathcal{S}}^{\phantom{i}} \YJ 
\mkern1mu h_{\mkern1mu i \mkern-1mu j} = \sum_{i \mkern-1mu j \mkern2mu 
\in \mathcal{S}} \sum_{k \mkern0.5mu l \mkern1mu \in \mathcal{S}} (c_i - 
c_k) (\CJ - c_l) \mkern1mu Z_{21 - i - k} \mkern1mu h_{\mkern1mu i 
\mkern-1mu j} \mkern1mu h_{kl} \\
   a_{\mkern0.5mu i \mkern1mu 4} &= \frac{1}{c_4^2} \biggl( c_i^2 (c_i 
- c_4) - 3 \sum_{j = 5}^{i - 1} \AIJ c_j (c_j - c_3) \biggr) , \qquad i 
\ge 4 \\
   a_{\mkern0.5mu i \mkern1mu 3} &= \frac{1}{c_3^2} \biggl( c_i^2 \bigl( 
c_4 - {\TS\frac{2}{3}} c_i \bigr) + 2 \sum_{j = 5}^{i - 1} \AIJ c_j (c_j 
- c_4) \biggr) , \qquad i \ge 4 \hskip24pt 
\raisebox{-13pt}{\smash{\begin{minipage}{0.15in}$\displaystyle\left. 
\vphantom{\rule{1pt}{60pt}} \right\}$\end{minipage}}} \mkern9mu 
\mbox{\small\smash{\begin{minipage}{1.26in} $\phantom{.}$ \vskip24pt 
ensures $\vphantom{|_|}$\\ $\vphantom{.}$\quad $\VEC{q}_0 = \VEC{0} 
\STRUT$\\ $\vphantom{.}$\quad $(\VEC{q}_1)_i = 0 \STRUT$ for all $i \ne 
2$\\ $\vphantom{.}$\quad $(\MAT{A} \mkern1mu \VEC{q}_1 \mkern-1.5mu )_i = 
0 \STRUT$ for all $i \ne 3$ \\ $\vphantom{.}$\quad $(\VEC{q}_2)_i = 0 
\STRUT$ for all $i \ne 2, 3$ \end{minipage}}} \\
   a_{\mkern0.5mu i \mkern1mu 2} &= \left\{ \begin{array}{ll} c_3^2 / 
2 c_2, & \quad i = 3 \\ 0, & \quad i \ne 3 \end{array} \right. \\
   a_{\mkern0.5mu i \mkern1mu 1} &= c_i - \sum_{j = 2}^{i - 1} \AIJ 
 \qquad \mbox{for all }i \end{align*} 
where $\mathcal{S} = \{ 65, 75, 85, 76, 86, 87 \}$. In particular, 
$a_{41} = c_4 / 4$ and $a_{43} = 3 c_4 / 4$.

The vectors $\VEC{q}_2$ and $\VEC{q}_2 \mkern1mu \VEC{c}$ are linear 
combinations of $\VEC{q}_1$ and $\MAT{A} \mkern1mu \VEC{q}_1$, also 
$(\MAT{A} \mkern1mu \VEC{q}_1 \mkern-1.5mu ) \mkern1mu \VEC{c} = c_3 
\mkern1.5mu \MAT{A} \mkern1mu \VEC{q}_1$. The four vectors $\VEC{q}_1$, 
$\MAT{A} \mkern1mu \VEC{q}_1$, $\MAT{A}^2 \VEC{q}_1$, and $\VEC{q}_3$ 
should be in the null space of the matrix $\MAT{B} = \bigl[ \SKJ \bigr]$, 
see Table~\ref{order_conditions}. The interpolant matrix $\MAT{B}$ is 
generated as
 \begin{gather}
   \MAT{B} = \left[ \begin{array}{ccccccccc}
     \,1\, & \,0\, & \,0\, & \,0\, & \,0\, & \,0\, & \,0\, & \,0\, & \,0\, \\
     0 & \frac{1}{2} & 0 & 0 & 0 & 0 & 0 & 0 & 0 \\
     0 & 0 & \frac{1}{3} & 0 & 0 & 0 & 0 & 0 & 0 \\
     0 & 0 & 0 & \frac{1}{4} & 0 & 0 & 0 & 0 & 0 \\
     0 & 0 & 0 & 0 & \frac{1}{5} & 0 & 0 & 0 & 0
   \end{array} \right]
 \underbrace{\left[ \begin{array}{ccccccccc}
     \VEC{1} & \VEC{c} & \VEC{c}^2 & \VEC{c}^3 & \VEC{c}^4 & \VEC{q}_1 & 
\MAT{A} \mkern1mu \VEC{q}_1 & \MAT{A}^2 \VEC{q}_1 & \VEC{q}_3
   \end{array} \right]}_{9 \times 9 \textrm{~matrix}}
   {{\vphantom{\big|}}^{-1}}
 \label{mat_B} \end{gather} Here the last components of the four vectors 
$\VEC{q}_1$, $\MAT{A} \mkern1mu \VEC{q}_1$, $\MAT{A}^2 \VEC{q}_1$, and 
$\VEC{q}_3$ are set to $0$, which is compatible with the order 
conditions. The interpolant $\VEC{\beta}(\theta) = \bigl[ \, \theta ~~~ 
\theta^2 ~~~ \theta^3 ~~~ \theta^4 ~~~ \theta^5 \, \bigr] \, \MAT{B}$ is 
continuously differentiable if
 \begin{gather*}
   \left[ \begin{array}{ccccc}
     1\; & \;0\; & \;0\; & \;0\; & \;0 \\
     1\; & \;1\; & \;1\; & \;1\; & \;1 \\
     1\; & \;2\; & \;3\; & \;4\; & \;5
   \end{array} \right] \MAT{B} =
   \left[ \begin{array}{ccccccccc}
     1\; & \;0\; & \;0\; & \;0\; & \;0\; & \;0\; & \;0\; & \;0\; & \;0 \\ 
     {\mkern1mu}b{\mkern-1mu}_1 & {\mkern1mu}b_2{\mkern-1mu} & 
{\mkern1mu}b_3{\mkern-1mu} & {\mkern1mu}b_4{\mkern-1mu} & 
{\mkern1mu}b_5{\mkern-1mu} & {\mkern1mu}b_6{\mkern-1mu} & 
{\mkern1mu}b_7{\mkern-1mu} & {\mkern1mu}b_8{\mkern-1mu} & \;0 \\
     0\; & \;0\; & \;0\; & \;0\; & \;0\; & \;0\; & \;0\; & \;0\; & \;1
   \end{array} \right]
 \end{gather*} The \FIRST{} and \THIRD{} rows of this equation are 
satisfied automatically due to how the first and last rows of the $9 
\times 9$ matrix in eq.~(\ref{mat_B}) look like. The \SECOND{} row is 
used to determine $\VEC{a}_{9*} = \VEC{b}$. The interpolant functions 
$\beta_2(\theta) \equiv 0$ and $\beta_3(\theta) \equiv 0$ as 
$(\VEC{q}_1)_i = 0$ and $(\MAT{A} \mkern1mu \VEC{q}_1 \mkern-1.5mu )_i = 
0$ for all $i \ne 2$ and $i \ne 3$, respectively.

The $5$-dimensional (could be parameterized by $c_2$, $c_4$, $c_5$, 
$c_6$, and $c_7$) family of $(5, 6)$ pairs (see \citep[][tab.~5]{DLMP89}, 
\citep[][tab.~2]{Ver91}, \citep[][tab.~4]{ShVe94}, 
\citep[][tab.~3]{Ver10}) is constructed in a similar fashion. It has 
$\textrm{DSO} = 3$, $c_8 = 1$, and $a_{i 2} = 0$ for all $i \ne 3$. The 
$9 \times 9$ matrix in the equation (\ref{mat_B}) is singular. With just 
$9$ stages, in order for the {\NTH{5}} order interpolant to exist, the 
last $4$ columns of the matrix should be linearly dependent, which 
happens if $c_4 = c_5$. Then the {\NTH{5}} stage repeats the {\NTH{4}} 
one, while some of the coefficients become infinite. To equip such an 
embedded $(5, 6)$ pair with an interpolant of order $5$, one needs to add 
at least one more stage.

\section{Choice of the degrees of freedom} \label{parameters}

A measure of the amount of non-positivity present in an interpolant is 
its total variation: \begin{gather*} V \vphantom{|}_0^1(\MAT{B}) = 
\sum_{j = 1}^s \, \int\limits_0^1 \d \theta \; \biggl| \frac{\d 
\SJ(\theta)}{\d \theta} \biggr| \end{gather*} As $\sum_j \SJ(\theta) 
\equiv \theta$, the minimal possible value of $V 
\vphantom{|}_0^1(\MAT{B})$ is equal to $1$. If the integration over 
$\theta$ is restricted to the region where $\SJ(\theta)$ has a negative 
slope: $N \vphantom{|}_0^1(\MAT{B}) = \sum_j \int_0^1 \d \theta \; \bigl| 
\beta'_{\mkern-0.25mu j}(\theta) \bigr| \, H \bigl( 
-\beta'_{\mkern-0.25mu j}(\theta) \bigr)$, then $V 
\vphantom{|}_0^1(\MAT{B}) = 1 + 2 N \vphantom{|}_0^1(\MAT{B})$. Here $H$ 
is the Heaviside step function.

\begin{figure} \centerline{\includegraphics{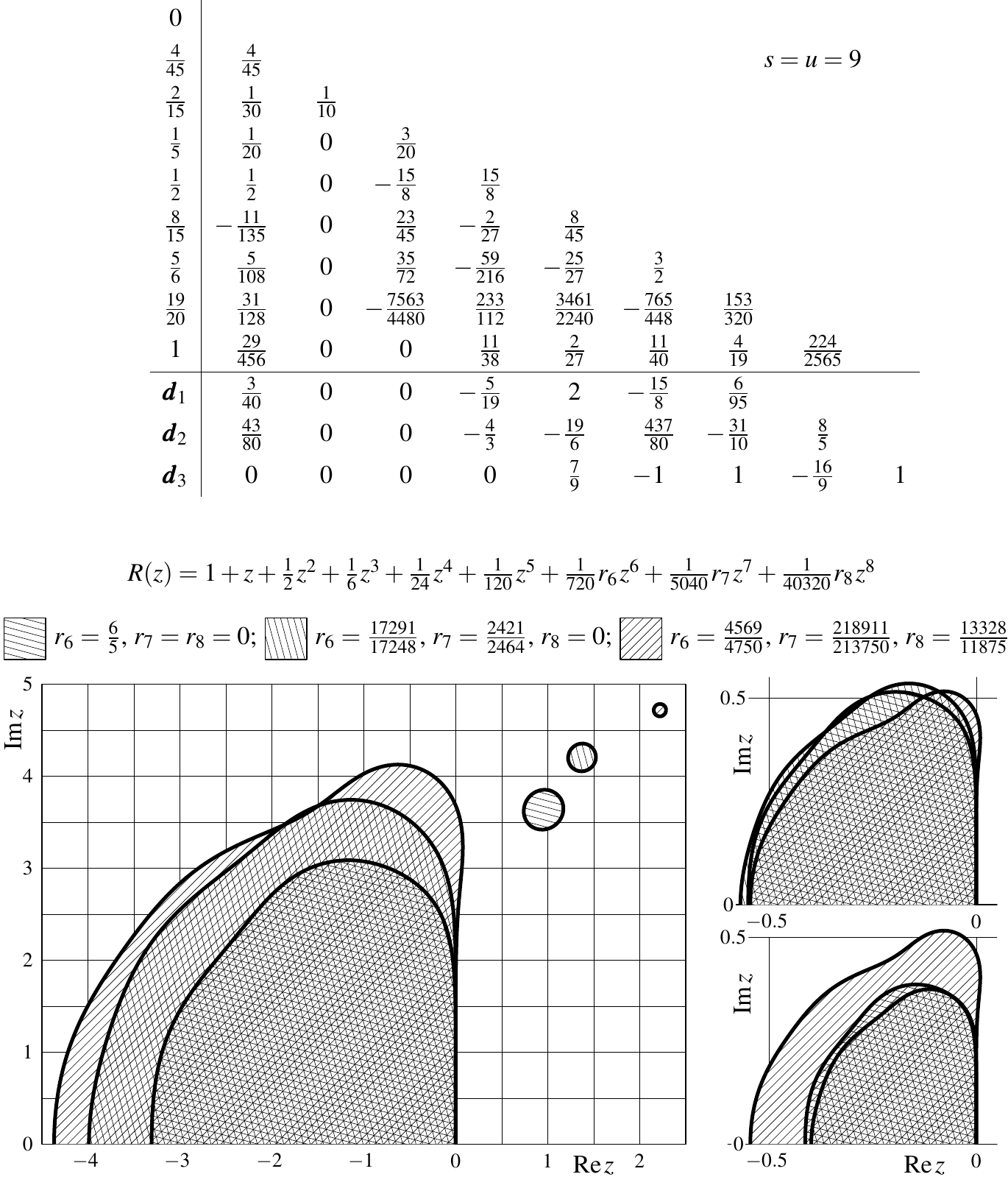}} 
\caption{A continuous $(4, 5)$ pair: the Butcher tableau and the region 
of absolute stability $\bigl| \mkern1mu R(z) \bigr| \le 1$. The weights 
vector $\VEC{b}$ repeats the last row of $\MAT{A}$ and is not shown. The 
difference between the weights vectors of the {\NTH{4}} and the {\NTH{5}} 
order methods within the pair could be any linear combination of 
$\VEC{d}_1$, $\VEC{d}_2$, and $\VEC{d}_3$. Regions of absolute stability 
of the \citep[][tab.~2]{DoPr80} pair (with $r_7 = 0$) and the 
\citep[][]{BoSh96} pair (with $r_7 = \frac{2421}{2464}$) are shown for 
comparison. The three regions of absolute stability shown on the small 
right panels are scaled for equal cost, \IE, the regions where $\bigl| 
\mkern1mu R \bigl( \mkern-1mu (s - 1) \mkern1mu z \bigr) \bigr| \le 1$ 
are depicted. For the \citep[][tab.~2]{DoPr80}, \citep[][]{BoSh96} pairs 
and the pair presented in this figure the number of stages $s$ is equal 
to $7$, $8$, $9$ on upper panel (just position is updated, $\VEC{x}(t) 
\mapsto \VEC{x}(t + h)$) and $9$, $11$, $9$ on lower panel (interpolant 
of order $5$ with $\max_{\, 0 \le \theta \le 1} T_6(\theta) \approx T_6$ 
is computed), respectively.} \label{pair1} \end{figure}
 \begin{figure} \centerline{\includegraphics{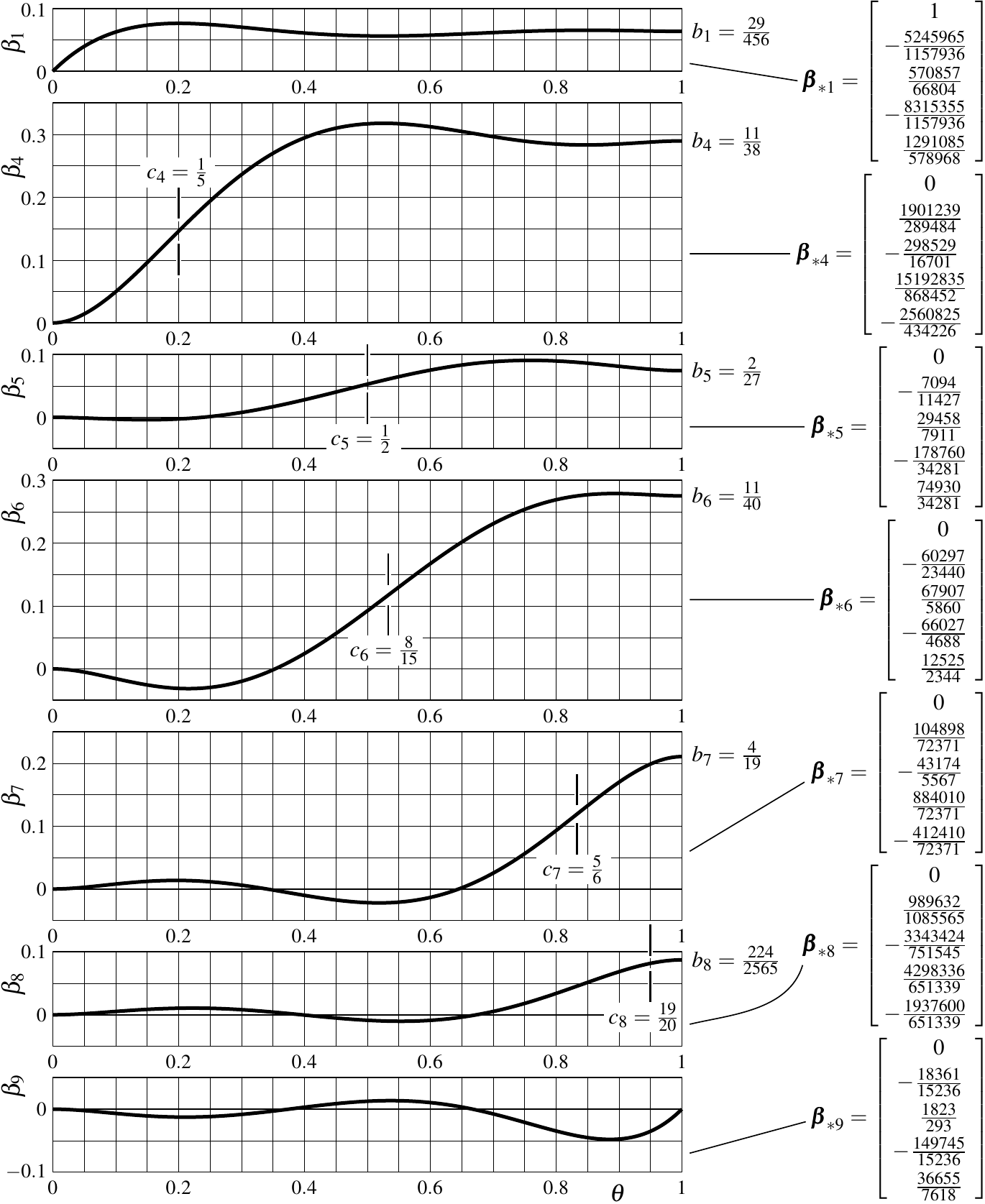}} 
\caption{Columns of the matrix $\MAT{B}$ and interpolant functions 
$\SJ(\theta)$, where $j = 1$, $4$, $5$, \DOTS, $9$, for the pair shown in 
Figure~\ref{pair1}.} \label{spline1} \end{figure}
 \begin{figure} \centerline{\includegraphics{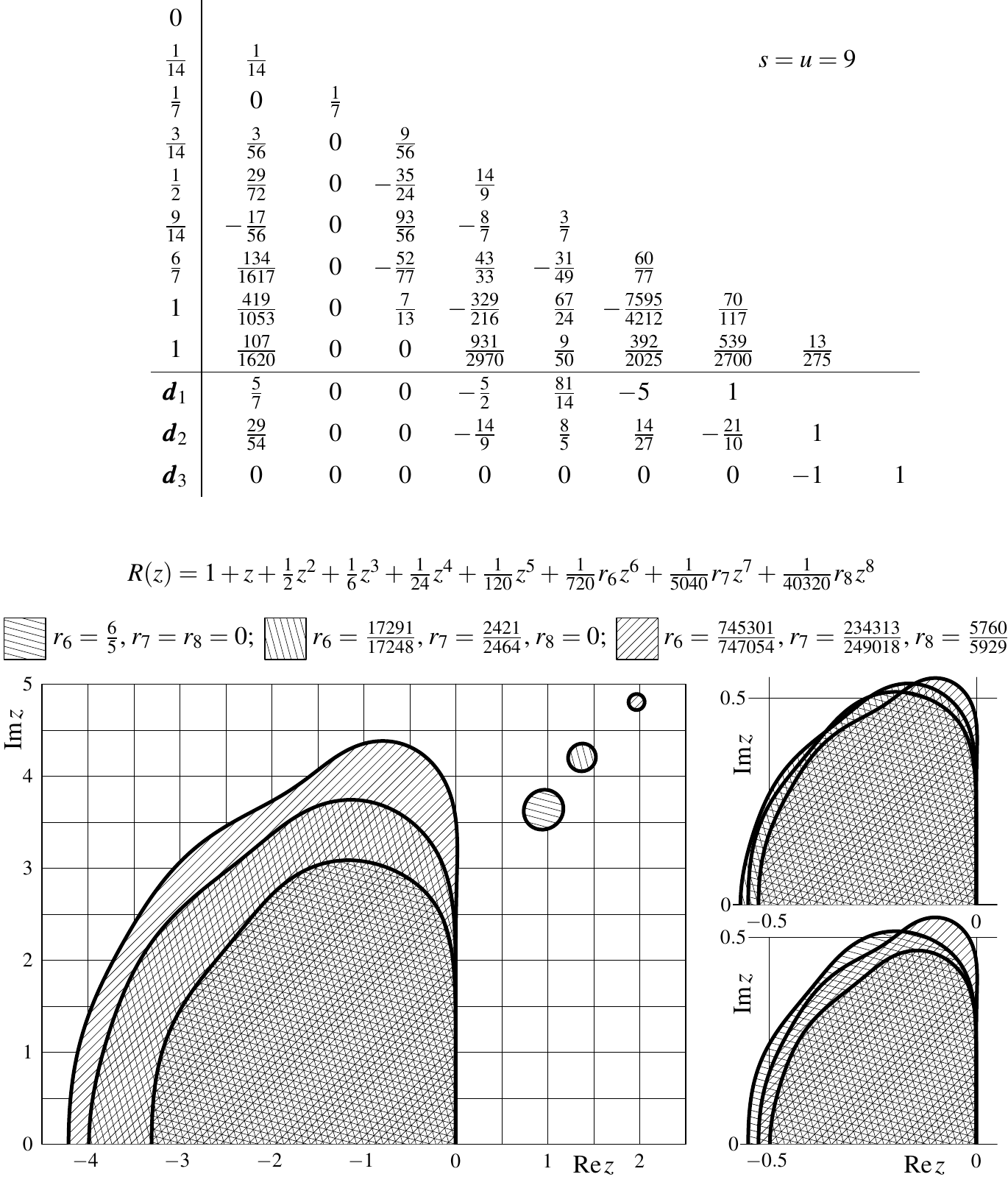}} 
\caption{A continuous $(4, 5)$ pair: the Butcher tableau and the region 
of absolute stability $\bigl| \mkern1mu R(z) \bigr| \le 1$. The weights 
vector $\VEC{b}$ repeats the last row of $\MAT{A}$ and is not shown. The 
difference between the weights vectors of the {\NTH{4}} and the {\NTH{5}} 
order methods within the pair could be any linear combination of 
$\VEC{d}_1$, $\VEC{d}_2$, and $\VEC{d}_3$. Regions of absolute stability 
of the \citep[][tab.~2]{DoPr80} pair (with $r_7 = 0$) and the 
\citep[][]{BoSh96} pair (with $r_7 = \frac{2421}{2464}$) are shown for 
comparison. The three regions of absolute stability shown on the small 
right panels are scaled for equal cost, \IE, the regions where $\bigl| 
\mkern1mu R \bigl( \mkern-1mu (s - 1) \mkern1mu z \bigr) \bigr| \le 1$ 
are depicted. For the \citep[][tab.~2]{DoPr80}, \citep[][]{BoSh96} pairs 
and the pair presented in this figure the number of stages $s$ is equal 
to $7$, $8$, $9$ on upper panel (just position is updated, $\VEC{x}(t) 
\mapsto \VEC{x}(t + h)$, or interpolant of order $4$, $4$, $5$ is 
computed) and $7$, $9$, $9$ on lower panel (low-cost interpolant of order 
$4$, $5$, $5$ is computed), respectively.} \label{pair2} \end{figure}
 \begin{figure} \centerline{\includegraphics{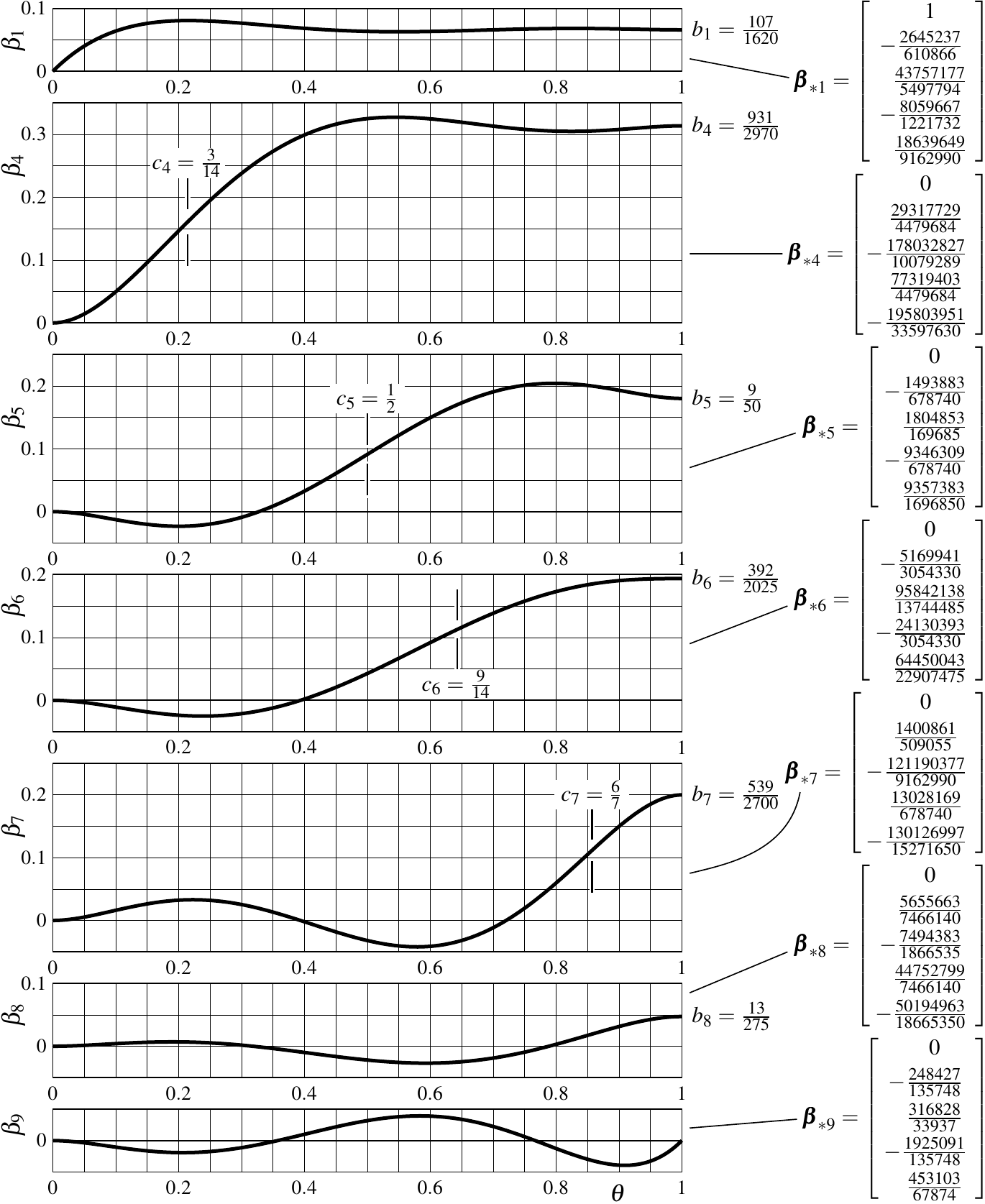}} 
\caption{Columns of the matrix $\MAT{B}$ and interpolant functions 
$\SJ(\theta)$, where $j = 1$, $4$, $5$, \DOTS, $9$, for the pair shown in 
Figure~\ref{pair2}.} \label{spline2} \end{figure}%
\hskip\parindent The continuous $(4, 5)$ pair shown in Figure~\ref{pair1} 
was obtained by minimizing the following function: \begin{gather*} 
\max_{0 \le \theta \le 1} T_6(\theta) + 10^{-4} \mkern2mu V 
\vphantom{|}_0^1(\MAT{B}) + 10^{-7} \sum_{ij} \bigl( 4 | \AIJ
|^2 + | \AIJ |^4 \bigr) \end{gather*} The term $\max_{0 \le \theta \le 1} 
T_6(\theta)$ was motivated by the discussion in \citep[][p.~24]{BoSh96}. 
In the resulted pair the local error $T_6$ is only slightly smaller than 
$\max_{\, 0 \le \theta \le 1} T_6(\theta)$, see 
Table~\ref{table_methods}. The term with $V \vphantom{|}_0^1(\MAT{B})$ 
makes the interpolant functions to wiggle less. For the pair to have 
rational coefficients, the values of the $11$ parameters obtained in the 
minimization were approximated by rational numbers.

The pair shown in Figure~\ref{pair2} was obtained by minimizing $T_6$ 
subject to the $T_7 \le 10 T_6$ inequality constraint. (The 
\citep[][tab.~2]{DoPr80} and \citep[][p.~20]{BoSh96} pairs have the ratio 
$T_7 / T_6$ close to $10$. For the latter pair this was a part of the 
design process.) A minor variation of the pair brings the local error 
$T_6$ to zero, see Table~\ref{6th_order}. The value of $\max_{0 \le 
\theta \le 1} T_6(\theta)$ is much larger than $T_6$, see 
Table~\ref{table_methods}. This by itself is not necessarily problematic, 
as less accurate interpolant is connecting endpoints of the integration 
step that contain an error accumulated in many steps. With $(4, 5)$ pairs 
it is not even uncommon to use interpolants of order $4$.

\begin{table} \begin{gather*} \begin{array}{r|lll|ll}
     & 10^5 \times T_6 & 10^5 \times T_7 & 10^5 \times \max_{\theta} 
T_6(\theta) & \max_{ij} | \mkern1mu \AIJ \mkern1mu | & \phantom{0}V 
\vphantom{|}_0^1(\MAT{B}) \STRUT \\
     \hline
     \mbox{Owren--Zennaro} & 108.62... & 154.05... & 
\phantom{00}108.62... & \phantom{0}3.75 & 1.6496... \STRUT \\
     \mbox{Dormand--Prince} & \phantom{0}39.908... & 
395.57... & \phantom{000}39.908... & 11.595... & 3.2490... \STRUT \\ 
     \mbox{Bogacki--Shampine} & {} \phantom{00}2.2169... & {} 
\phantom{0}21.260... & \phantom{0000}2.2169... & \phantom{0}1.1637... & 
3.3092... \STRUT \\
     \hline \mbox{Figure~}\ref{pair1} & \phantom{00}9.2847... & 
\phantom{0}19.904... & \phantom{0000}9.7178... & \phantom{0}2.0803... & 
1.4857... \STRUT \\
     \mbox{Figure~}\ref{pair2} & \phantom{00}0.59809... & 
\phantom{00}5.9203... & \phantom{000}16.134... & \phantom{0}2.7916... & 
1.6424... \STRUT
   \end{array} \\[-28pt] \phantom{.} \end{gather*} \caption{A comparison 
of five continuous $(4, 5)$ pairs. The first three are from the 
literature: \citep[][fig.~3]{OwZe92}, \citep[][tab.~2]{DoPr80}, and 
\citep[][]{BoSh96}.} \label{table_methods} \end{table}

\begin{table} \begin{gather*} \begin{array}{r|lll}
     & 10^5 \times T_5(\VEC{x}) & 10^5 \times T_6(\VEC{x}) & 10^5 \times 
T_7(\VEC{x}) \STRUT \\
     \hline
     \mbox{Dormand--Prince, } \VEC{x} = 
\VEC{b}_{4^{\textrm{th}}\textrm{~order}} & 118.29... & 182.37... & 
414.05...\STRUT \\
     \VEC{x} = \frac{1}{3} \VEC{b} + \frac{2}{3} 
\VEC{b}_{4^{\textrm{th}}\textrm{~order}} & \phantom{0}78.863... & 
118.66... & 392.39...\STRUT \\
     \hline
     \mbox{Bogacki--Shampine, }\VEC{x} = \VEC{b} + \VEC{E} & 
\phantom{0}10.595... & {} \phantom{0}12.204... & \phantom{0}24.114... 
\STRUT \\
     \VEC{x} = \VEC{B} & \phantom{0}10.615... & \phantom{0}10.992... & 
\phantom{0}20.562... \STRUT \\
     \hline
     \mbox{Figure~}\ref{pair1}\mbox{, }\VEC{x} = \VEC{b} + \frac{1}{42} 
\VEC{d}_1 & \phantom{0}19.765... & \phantom{0}14.406... & 
\phantom{0}18.948...\STRUT \\
     \VEC{x} = \VEC{b} - \frac{1}{323} \VEC{d}_2 & \phantom{0}19.465... 
& \phantom{0}29.174... & \phantom{0}33.473...\STRUT \\
     \VEC{x} = \VEC{b} + \frac{1}{152} \VEC{d}_3 & \phantom{0}19.511... 
& \phantom{0}27.858... & \phantom{0}27.751...\STRUT \\
     \hline
     \mbox{Figure~}\ref{pair2}\mbox{, }\VEC{x} = \VEC{b} + \frac{1}{2430} 
\VEC{d}_1 & \phantom{00}0.99808... & \phantom{00}1.0189... & 
\phantom{00}6.0611...\STRUT \\
     \VEC{x} = \VEC{b} + \frac{1}{2156} \VEC{d}_2 & 
\phantom{00}0.99189... & \phantom{00}1.4234... & 
\phantom{00}6.4042...\STRUT \\
     \VEC{x} = \VEC{b} - \frac{1}{986} \VEC{d}_3 & 
\phantom{00}0.99041... & \phantom{00}1.3645... & 
\phantom{00}6.4662...\STRUT
   \end{array} \\[-28pt] \phantom{.} \end{gather*} \caption{Local errors 
for the lower order methods of four embedded $(4, 5)$ pairs. The second 
line of the \citep[][tab.~2]{DoPr80} entry is the modification suggested 
in \citep[][p.~141]{Sha86}. The notation $\VEC{E}$ and $\VEC{B}$ in the 
\citep[][]{BoSh96} entry is taken from the 
\scalebox{0.85}[1]{\tt\bfseries rksuite.f} code \citep[][]{RKSUITE}.} 
\label{table_methods_d} \end{table}

\begin{table} \begin{gather*} \begin{array}{c|ccccccccc}
     0\STRUT \\
     \frac{1}{14}\STRUT & \phantom{-}\frac{1}{14} & & & & & & & 
\mathclap{s = u = 9} \\
     \frac{1}{7}\STRUT & \phantom{-}0 & \phantom{-}\frac{1}{7} \\
     \frac{3}{14}\STRUT & \phantom{-}\frac{3}{56} & \phantom{-}0 & 
\phantom{-}\frac{9}{56} \\
     \frac{1}{2}\STRUT & \phantom{-}\frac{29}{72} & \phantom{-}0 & 
-\frac{35}{24} & \phantom{-}\frac{14}{9} \\
     \frac{9}{14}\STRUT & -\frac{17}{56} & \phantom{-}0 & 
\phantom{-}\frac{93}{56} & -\frac{8}{7} & \phantom{-}\frac{3}{7} \\
     \frac{6}{7}\STRUT & \phantom{-}\frac{199}{1372} & \phantom{-}0 & 
-\frac{195}{196} & \phantom{-}\frac{1259}{784} & -\frac{3855}{5488} & 
\phantom{-}\frac{45}{56} \\
     1\STRUT & \phantom{-}\frac{4903}{25596} & \phantom{-}0 & 
\phantom{-}\frac{4487}{2844} & -\frac{255101}{102384} & 
\phantom{-}\frac{33847}{11376} & -\frac{94325}{51192} & 
\phantom{-}\frac{3773}{6399} \\
     1\STRUT & \phantom{-}\frac{16}{243} & \phantom{-}0 & \phantom{-}0 & 
\phantom{-}\frac{16807}{53460} & \phantom{-}\frac{53}{300} & 
\phantom{-}\frac{2401}{12150} & \phantom{-}\frac{2401}{12150} & 
\phantom{-}\frac{79}{1650} \\
 \end{array} \\[-28pt] \end{gather*} \caption{The Butcher tableau of the 
embedded $(4, 6)$ pair with FSAL property. The weights vector $\VEC{b}$, 
which produces the {\NTH{6}} order update, repeats the last row of 
$\MAT{A}$ and is not shown. The $9 \times 9$ matrix in eq.~(\ref{mat_B}) 
is non-singular, and no weights vector other than $\VEC{b}$ would give a 
method of order $5$. The {\NTH{5}} order interpolant is constructed as in 
eq.~(\ref{mat_B}). The differences $\VEC{d}_{1,2,3}$ between the weights 
vectors of the {\NTH{4}} and the {\NTH{6}} order methods, that could be 
used for error control, are the same as in Figure~\ref{pair2}.} 
\label{6th_order} \end{table}

\section{Numerical tests} \label{tests}

\begin{figure} \centerline{\includegraphics{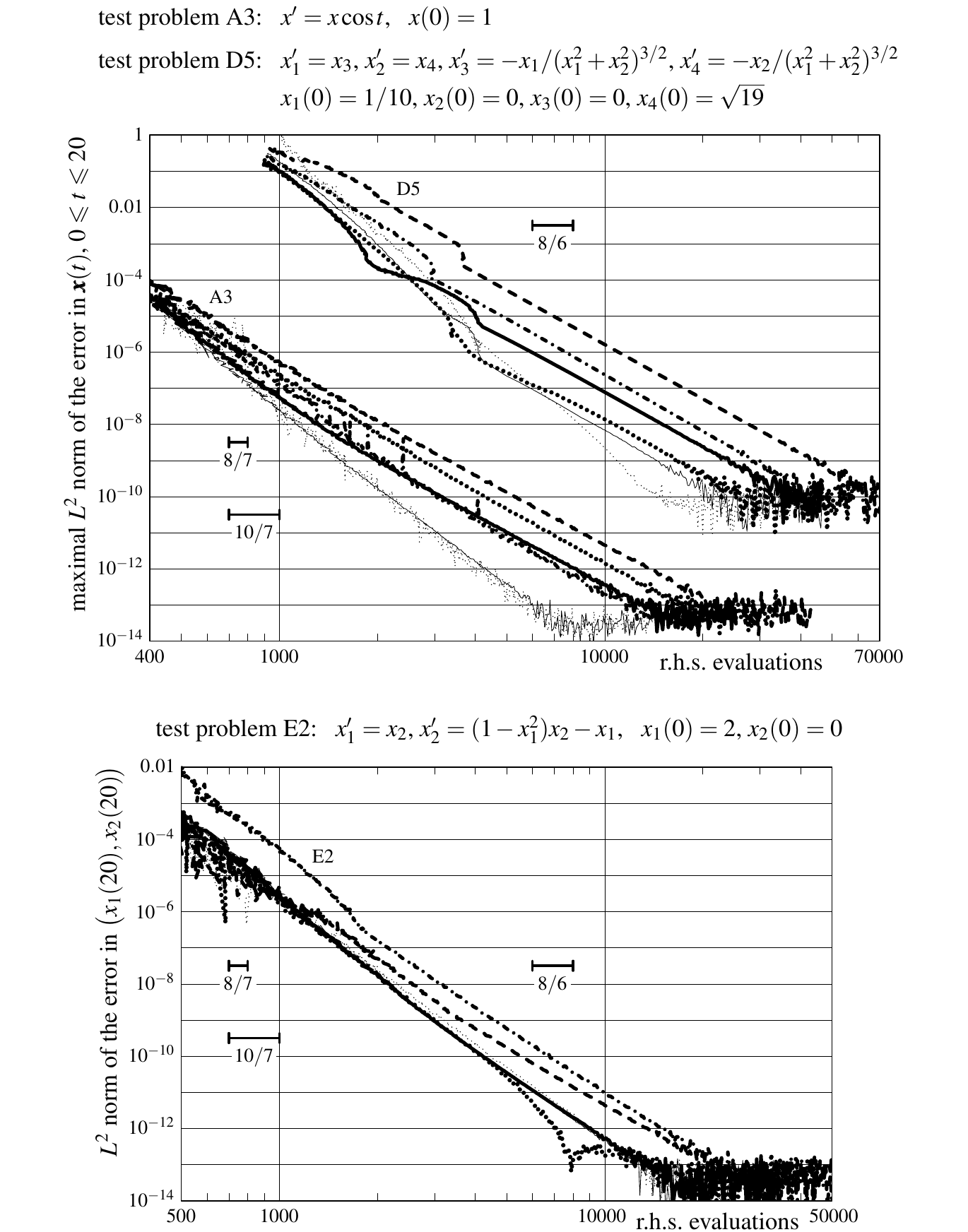}} 
\caption{Efficiency curves for problems~A3 \citep[][p.~617]{DETEST}, D5 
\citep[][p.~620]{DETEST}, and E2 \citep[][p.~621]{DETEST}: the 
\citep[][tab.~2]{DoPr80} pair (dash-dotted curve), the \citep[][]{BoSh96} 
pair (dotted curve), the pair in Figure~\ref{pair1} (dashed curve), the 
pair in Figure~\ref{pair2} (solid curve), the $(4, 6)$ pair in 
Table~\ref{6th_order} (thin solid curve), and the \citep[][tab.~3]{Ver10} 
par (thin dotted curve).} \label{A3_D5_E2} \end{figure}

\begin{figure} \centerline{\includegraphics{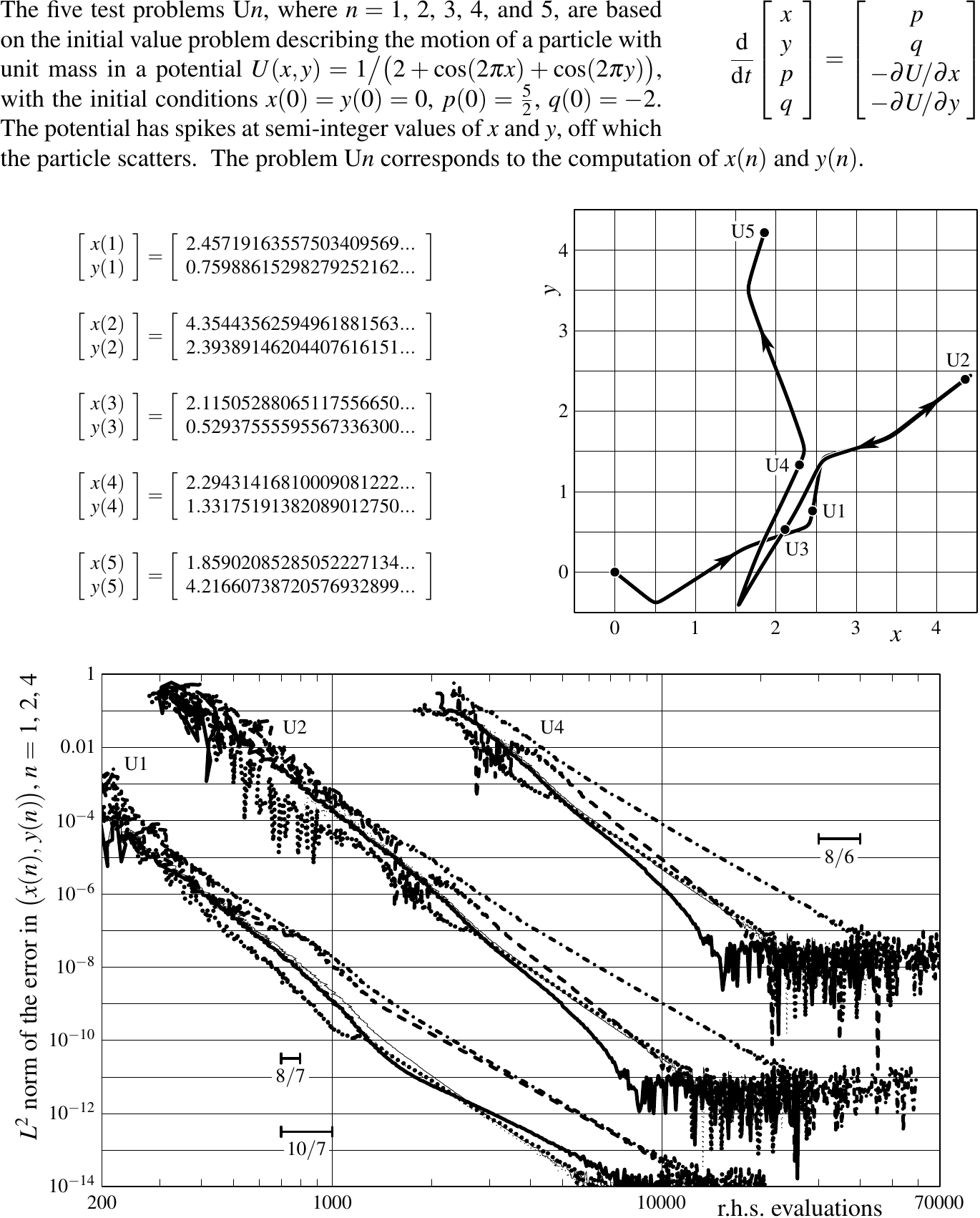}} 
\caption{The test problems U$n$, where $n = 1$, $2$, $3$, $4$, and $5$; 
accurate up to $10^{-20}$ values of $\smash{\bigl( x(n), \, y(n) 
\bigr)}$, $n = 1$, $2$, $3$, $4$, and $5$ (middle left); the trajectory 
$\bigl( x(t), \mkern1.5mu y(t) \mkern-1mu \bigr)$ for $0 \le t \le 5$ 
(middle right); and efficiency curves for problems U1, U2, and U4 (bottom 
panel): the \citep[][tab.~2]{DoPr80} pair (dash-dotted curve), the 
\citep[][]{BoSh96} pair (dotted curve), the pair in Figure~\ref{pair1} 
(dashed curve), the pair in Figure~\ref{pair2} (solid curve), the $(4, 
6)$ pair in Table~\ref{6th_order} (thin solid curve), and the 
\citep[][tab.~3]{Ver10} par (thin dotted curve).} \label{U1_U2_U4} 
\end{figure}

The performance of the pairs constructed in the previous section is 
demonstrated on test problems A3, D5, and E2 from \citep[][]{DETEST} in 
Figure~\ref{A3_D5_E2}; and on new suggested test problems U1, U2, and U4 
in Figure~\ref{U1_U2_U4}. The error, that is the difference between the 
exact and numerical solutions, was computed only at the ends of 
integration steps. The adaptive step size scheme $h \leftarrow 0.9 
\mkern1mu h \mkern1.5mu (\mbox{ATOL} / E)^{1 / 5}$ was used. (The 
starting step size $h_0 = 10^{-3}$ was swiftly corrected by the adaptive 
step size control.) Here $\mbox{ATOL}$ is the absolute error tolerance, 
and $E$ is the $l^2$-norm of the difference vector between the two 
solutions within a pair. The steps with $E > \mbox{ATOL}$ were rejected, 
but they were still contributing to the number of the r.h.s.~evaluations. 
For the pairs with multiple difference vectors between the weights of the 
higher and lower order methods (\IE, the \citep[][]{BoSh96} pair and the 
pairs in Figures~\ref{pair1} and \ref{pair2}, see 
Table~\ref{table_methods_d}), the difference vectors that use smaller 
number of stages were tried first. Once a step was rejected, no further 
difference vectors were tried. The size of the next step, whether the 
previous step was rejected or not, was chosen according to the maximal 
value of $E$ between the tried vectors.

The pairs with higher order $6$ performed better on test problems A3 and 
D5. On other problems their performance was similar to the 
\citep[][]{BoSh96} pair. The \citep[][tab.~2]{DoPr80} pair performed well 
on problem A3, while on problems E2, U2, and U4 it performed the worst. 
The pair shown in Figure~\ref{pair1} performed the worst on problems A3 
and D5. The pair from Figure~\ref{pair2} seems to be at least as 
efficient as the \citep[][]{BoSh96} pair. Note that the efficiency curves 
in Figures~\ref{A3_D5_E2} and \ref{U1_U2_U4} show the cost of obtaining 
the numerical solution without computing interpolants. In case of the 
\citep[][]{BoSh96} pair, no additional stages are needed to use the 
interpolant of order $4$ \citep[][]{Bog90}. One or three additional 
r.h.s.~evaluations per step are needed to compute less and more accurate 
interpolant of order $5$, respectively \citep[][p.~24]{BoSh96}, which 
corresponds to the increase of the cost by factors $8 / 7$ and $10 / 7$.

 \begin{figure} \centerline{\includegraphics{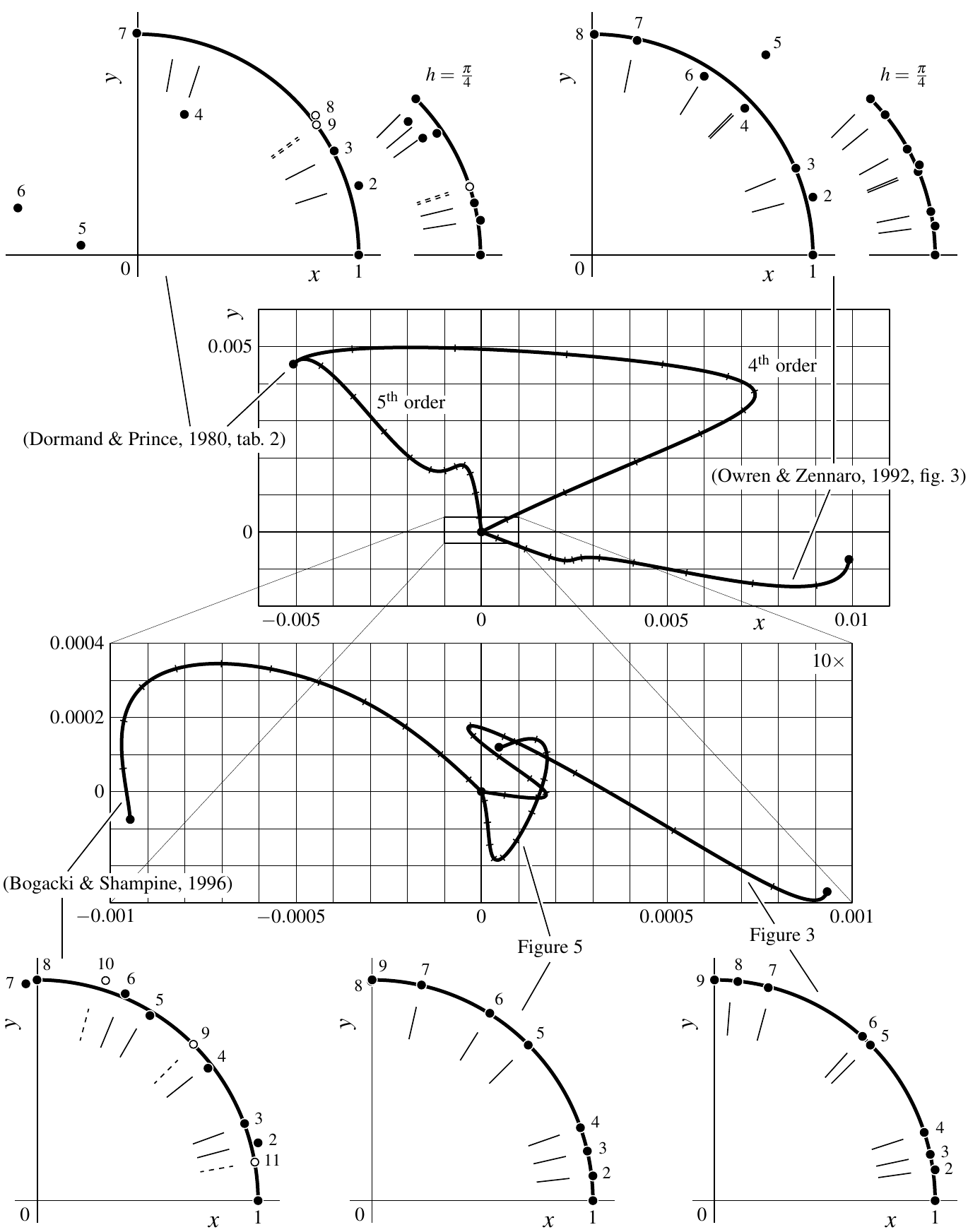}} 
\caption{One step $h = \frac{\pi}2$ for the system $\d x / \d \mkern1mu t 
= -y$, $\d y / \d \mkern1mu t = x$ with initial conditions $x(0) = 1$, 
$y(0) = 0$: the intermediate positions $\VEC{X}_{\mkern-1mu i}$, $1 \le i 
\le s$, and the error $\bigl( x(\frac{\pi}{2} \theta) - 
\cos(\frac{\pi}{2} \theta), \, y(\frac{\pi}{2} \theta) - 
\sin(\frac{\pi}{2} \theta) \bigr)$, where $0 \le \theta \le 1$, made by 
an interpolant. Nodes are shown by the radial ticks. Open circles 
correspond to additional stages that are used to construct an 
interpolant. Ticks on the error curves correspond to $\theta = 
\frac1{12}$, $\frac16$, $\frac14$, \DOTS, $\frac{11}{12}$. In case of the 
\citep[][tab.~2]{DoPr80} pair, the \NTH{5} order interpolant is the 
\citep[][]{CMR90} one, while the \NTH{4} order interpolant with $u = s = 
7$ is from \scalebox{0.85}[1]{\tt\bfseries ntrp45.m} code that is a part 
of MATLAB\textsuperscript{\textregistered} software.} \label{circle} 
\end{figure}%
\hskip\parindent{\!}Figure~\ref{circle} shows the performance of five 
pairs and their interpolants on the system $\d x / \d t = -y$, $\d y / \d 
t = x$ with initial condition $\bigl( x(0), y(0) \bigr) = (1, 0)$. Just 
one step $h = \frac{\pi}2$ is made, with no error control. The exact 
solution is $x(t) = \cos t$, $y(t) = \sin t$ with $\bigl( x(\frac{\pi}2), 
y(\frac{\pi}2) \bigr) = (0, 1)$. In cases of the \citep[][tab.~2]{DoPr80} 
and \cite[][fig.~3]{OwZe92} pairs the intermediate positions 
$\VEC{X}_{\mkern-1mu i}$, $1 \le i \le s$, are also shown for the step 
size $h = \frac{\pi}4$. For the \citep[][tab.~2]{DoPr80} pair, 
$(\VEC{q}_2)_i = 0$ for all $i \ne 2$. In the leading order the deviation 
of $\VEC{X}_{\mkern-1mu i}$ from $\VEC{x}(t + c_{\mkern0.5mu i} \mkern1mu 
h)$ is controlled by $T_3(\VEC{a}_{i*}, c_i) = \bigl| \mkern0.5mu 
\VEC{a}_{i*} \mkern1mu \MAT{A} \mkern1mu \VEC{c} - c_i^3 / 6 \bigr|$ that 
for $i = 3$, $4$, $5$, and $6$ is equal to $\frac{9}{2000}= 0.0045$, 
$\frac{28}{375} = 0.0746...$, $\frac{2536}{10935} = 0.231...$, and 
$\frac{71}{330} = 0.215...$, respectively. Relatively large values of 
$T_3(\VEC{a}_{5*}, c_5)$ and $T_3(\VEC{a}_{6*}, c_6)$ explain the 
observed deviation of $\VEC{X}_{\mkern-1mu 5}$ and $\VEC{X}_{\mkern-1mu 
6}$ from $\bigl( \cos \frac{8 \mkern0.5mu h}{9}, \, \sin \frac{8 
\mkern0.5mu h}{9} \bigr)$ and $(\cos h, \, \sin h)$. Even for $h = 
\frac{\pi}2$ the position update $\VEC{x}(t + h) = \VEC{X}_{\mkern-1mu 
7}$ is remarkably close to the exact value, though.

\section{Conclusion}

Utilizing $9$ stages, it is possible to construct embedded $(4, 5)$ 
pairs of explicit Runge--Kutta methods with FSAL property that are as 
cost-efficient as the best known conventional (\IE, interpolant is 
either of order $4$ or would require extra stages) pairs (\EG, the 
Dormand--Prince and Bogacki--Shampine ones), but with the benefit of 
having continuous formulae or interpolants of order $5$ available at no 
additional cost.

 \bibliographystyle{IMANUM-BIB}

\end{document}